\documentclass[10pt]{amsart}
\usepackage{amsmath}
\usepackage{amssymb}
\usepackage{enumerate}
\usepackage{amsbsy}
\usepackage{amsfonts}
\usepackage{color}
\usepackage{upgreek}

\newtheorem{thm}{Theorem}[section]
\newtheorem{lem}[thm]{Lemma}
\newtheorem{prop}[thm]{Propsition}
\newtheorem{cor}[thm]{Corollary}
\newtheorem{defn}[thm]{Definition}
\newtheorem{rem}[thm]{Remark}

\numberwithin{equation}{section}

\begin{document}

	\title[Multilinear estimates with angular regularity]{Improved multilinear estimates and global regularity for general nonlinear wave equations in $(1+3)$ dimensions}

    \author[S. Hong]{Seokchang Hong}
    \address{Department of Mathematics, Chung-Ang University, Seoul 06974, Korea}
    \email{seokchangh11@cau.ac.kr}

	
	\begin{abstract}
	This paper is devoted to the investigation of long-time behaviour of solutions to wave equations with quadratic nonlinearity and cubic Dirac equations with Hartree-type nonlinearity. 
We consider the nonlinearity here with enough simplicity so that we can treat it as a toy model and simultaneously with enough generality so that we can apply our result to wave and Dirac equations with various nonlinearities. 	
	The challenging point is that nonlinearity possesses singularity near the origin. Our strategy is to relax such a singularity by exploiting fully an angular momentum operator. In this manner we establish scattering for the critical Sobolev data. 
	%The nonlinearity here contains enough simplicity so that we treat it as a toy model and at the same time it contains enough generality so that we can apply our result to wave and Dirac equations with various nonlinearities. 	
	\end{abstract}

		\maketitle

\section{Introduction}
For several decades wave equations have appealed to a lot of interest and have been extensively studied in many works of literature. It plays a role as a model to explain various physical phenomena and in the mathematical literature, the study of the wave equations becomes the very first step to shattering the light on the investigation of hyperbolic partial differential equations. 

We are interested in time-evolution of solutions to wave equations with various nonlinearities for low regularity initial data. In the investigation it is important to control the nonlinearity in terms of the inital data. In other words, we have to prove that the presence of nonlinearity turns out to be nothing but a small perturbation. Such a perturbative method to wave equations and even more general dispersive equations is a typical approach to the study of Cauchy problems. The first well-known tool is so-called Stricharz estimates \cite{keeltao,stri}
$$
\|e^{-it|\nabla|}P_1f\|_{L^q_tL^r_x(\mathbb R^{1+n})} \lesssim \|P_1f\|_{L^2_x}
$$
for any function $f$. Here $P_1f$ is the projection onto the unit frequency and $(q,r)$ is a proper admissible pair. However, the linear estimate is not sufficient to control all over the frequency-interactions between the products of homogeneous solutions especially when we are concerned with well-posedness problem for a low regularity data. This problem requires one to delicately consider the following bilinear estimates
$$
\|e^{-it|\nabla|}f e^{-it|\nabla|}g\|_{L^2_tL^2_x} \lesssim \|f\|_{L^2_x}\|g\|_{L^2_x}.
$$
In fact, when nonlinearity is given by power-type, nonlinear estimates are reduced to bilinear estimates. Recently there has been a huge amount of progress on bilinear estimates of wave-type by many works of \cite{fosklai,klaima,klaima1,sleevar,tao,tao1,tao2,tataru,tataru1} and long-time behaviour of solutions to nonlinear wave equations and even more complicated systems such as the Maxwell-Klein-Gordon or the Yang-Mills equations is well-known for $(1+4)$ dimensions and higher dimensional setting \cite{kriegersterbenztataru, kriegertataru,sterbenz2,ohtataru}. 

However, global regularity is still open for the most of wave equations in a low dimensional setting such as $(1+3)$ or $(1+2)$ dimensions. At a first glimpse this is obviously because of the weaker time-decay of solutions in a low dimensional setting $\|e^{-it|\nabla|}P_1f\|_{L^\infty_x(\mathbb R^n)}\lesssim t^{-\frac{n-1}{2}}\|P_1f\|_{L^1_x}$. Moreover, at the nonlinear level, one can see that {\it resonant interactions} grows stronger as the spatial dimensions decrease. Even further, when nonlinearity possesses singularity near the origin, one may encounter more serious situation since the singularity grows harsher in a low dimension. From these several problems one may have a question whether it is possible to establish global well-posedness and scattering for the scale-invariant Sobolev data. 
To overcome this difficulty we equip the Sobolev spaces with an extra weighted smoothness assumptions with respect to the angular variables. Indeed, we invoke the infinitesimal rotation generators $\Omega_{ij}=x_i\partial_j-x_j\partial_i$. In the sprit of \cite{ster}, $\Omega_{ij}$ plays a crucial role in the aspect of both linear and multilinear estimates. More precisely, one enjoy a significant improvement of linear estimates. At the nonlinear estimates, the rotation operator helps to overcome the resonant interactions. Even more, the rotation can relax the harshness of the singularity since the operator $\Omega_{ij}$ works very favourably in the low-output interactions. In this manner, it is possible to improve the bilinear estimates.  

Now we turn to an application of an improved multilinear estimate on $\mathbb R^{1+3}$. We are concerned with somewhat a general class of quadratic nonlinear wave equations and the Hartree-type nonlinear Dirac equations which becomes {\it a toy model} for several nonlinear wave and Dirac equations. The following equations we shall present seem {\it too primitive} at a first glimpse, however, by the primitiveness and generality of a toy model we can attack efficiently even more complicated system such as gauge-field-theoretic wave equations which represent a genuinely physical model.
 %The application of the study on a toy model is the Maxwell-Klein-Gordon and the Yang-Mills equations in the Coulomb gauge and the Hartree-type nonlinear Dirac equations with more specific structures.
\subsection{Quadratic nonlinear wave equations}     
%The following two systems are our main concern throughout this paper and their nonlinearity possesses This paper is devoted to the study on the Cauchy problems of nonlinear wave equations. 
Firstly we aim to investigate global-in-time evolution of wave equations in $\mathbb R^{1+3}$ with quite a general quadratic nonlinearity given by
\begin{align}\label{main-wave}
\left\{
\begin{array}{l}
	\Box u = |\nabla|^{-1}Q(\overline u,u), \\
	(u,\partial_tu)|_{t=0} = (u_0,u_1),
\end{array}
\right.	
\end{align}
where $u$ is a complex-valued function on $\mathbb R^{1+3}$ and $Q:(u,v)\mapsto Q(u,v)$ is a bilinear form which is a finite linear combination of the standard $Q$-type null forms
\begin{align*}
Q_{ij}(u,v) & = \partial_iu\partial_jv - \partial_ju\partial_iv, \, Q_{0}(u,v)  = \partial_tu\partial_tv -\nabla u\cdot\nabla v,
\end{align*}
which give the cancellation by angle between input-frequency.\footnote{In fact, the null form $Q_0$ gives stronger cancellation, and we can overcome the singularity $|\nabla|^{-1}$ more easily by exploiting the $Q_0$ null form. However, for the generality of our result, we focus on the $Q_{ij}$ null form.} More precisely the Fourier transform of $Q_{ij}(u,v)$ is
\begin{align*}
|\widehat{Q_{ij}(u,v)}|(\zeta) \lesssim \int_{\zeta=\xi+\eta}\angle(\xi,\eta)|\xi||\eta| 	\widehat u(\xi) \widehat v(\eta)\,d\xi d\eta.
\end{align*}
The wave equation \eqref{main-wave} has the scaling symmetry, i.e., if $u=u(t,x)$, $(t,x)\in\mathbb R^{1+3}$ is a solution of \eqref{main-wave} then the scaled function $\lambda^{-1}u(\lambda^{-1}t,\lambda^{-1}x)$ will be also a solution to the equation \eqref{main-wave} for any $\lambda>0$ and hence the scale-invariant Sobolev space for the initial data $(u_0,u_1)$ is $\dot H^\frac12 \times \dot H^{-\frac12}$, where $\dot H^s$ is the usual homogeneous Sobolev space. Now we define the angularly regular space $\dot H^s_\sigma$ to be $\|f\|_{\dot H^s_\sigma}=\|\langle\Omega\rangle^\sigma f\|_{\dot H^s}$, where $\langle\Omega\rangle^\sigma=(1-\Delta_{\mathbb S^2})^\frac\sigma2$ and $\Delta_{\mathbb S^2}$ is the Laplace-Beltrami operator on the unit sphere $\mathbb S^2\subset\mathbb R^3$. The inhomogeneous Sobolev space with angular regularity $H^s_\sigma$ is defined in the obvious way. We state our first main result.
\begin{thm}\label{gwp-wave}
Let $\sigma=1$. Suppose that the initial datum $(u_0,u_1)\in \dot H^{\frac12}_\sigma\times\dot H^{-\frac12}_\sigma$ satisfies 
$$
\|(u_0,u_1)\|_{\dot H^\frac12_\sigma\times\dot H^{-\frac12}_\sigma} = \|u_0\|_{\dot H^\frac12_\sigma}+\|u_1\|_{\dot H^{-\frac12}_\sigma}\ll1.$$
 The Cauchy problem for the equation \eqref{main-wave} is globally well-posed and scatters to free solutions as $t\rightarrow\pm\infty$.
\end{thm}
\subsubsection{Application to the Maxwell-Klein-Gordon equations in the Coulomb gauge}
%Note that Theorem \ref{gwp-wave} means that flow map $u_0\mapsto u(t)$ is continuous in $C(\mathbb R; \dot H^\frac12_\sigma(\mathbb R^3;\mathbb C))$. More precisely we obtain solutions in the $U^2$ spaces. (See Section \ref{ftn-sp}.) The application of Theorem \ref{gwp-wave} we have in mind is the Maxwell-Klein-Gordon equations in the Coulomb gauge. 
%Here $u$ is a complex-valued function. 
We would like to mention here briefly an application of Theorem \ref{gwp-wave}. The Maxwell-Klein-Gordon system is a physical model for the interaction of a spin $0$ particle with electromagnetic fields. We define the real-valued gauge potentials $A_\mu$, $\mu=0,1,\cdots,3$ on the Minkowski space $(\mathbb R^{1+3},\mathbf m)$, where the metric $\mathbf m$ is given by $\mathbf m = \textrm{diag}(-1,1,1,1)$. The covariant derivative is given by $\mathcal D_\mu = \partial_\mu+iA_\mu $. The electromagnetic field $F$ associated to the potential $A_\mu$ is defined by $F_{\mu\nu}=\partial_\mu A_\nu-\partial_\nu A_\mu$. Then the covariant form of the Maxwell-Klein-Gordon system presents
\begin{align}
\begin{aligned}
	\partial_\mu F^{\mu\nu} &= -\textrm{Im}(\phi \overline{\mathcal D^\nu\phi}), \\
	\mathcal D_\mu\mathcal D^\mu\phi &= 0,
\end{aligned}	
\end{align}
where $\textrm{Im}(A)$ is the imaginary part of $A$. 
 Note that we adapt the usual summation convention with respect to repeated indices. The Maxwell-Klein-Gordon system has a gauge-invariance. Indeed, if $(A_\mu,\phi)$ is a solution to the system, then for any real-valued smooth function $\Lambda$ on $\mathbb R^{1+3}$, the set $(A_\mu+\partial_\mu\Lambda, e^{-i\Lambda}\phi)$ is also a solution to the system. This observation allows one to enjoy the gauge-freedom. Now we impose the Coulomb gauge condition: $\textrm{div} A = \partial_j A^j=0$. Then after an application of the projection $\mathbf P=-\frac{(\rm curl)^2}{\Delta}$ we see that the spatial parts of the gauge potentials obey the following wave equation
 \begin{align}\label{model-mkg}
 \Box A_j = -\textrm{Im}\,\mathbf P(\phi\overline{\mathcal D_j\phi}).
\end{align}
 Then the quadratic nonlinearity in the wave equation \eqref{model-mkg} presents a finite linear combination of the $Q$-type null forms as $\Delta^{-1}\partial_k Q_{ij}(\phi,\overline\phi)$, which turns out to be the nonlinearity in our toy model \eqref{main-wave}. We refer the readers to \cite{masterbenz} for more details on the Maxwell-Klein-Gordon system. 
 
 The Maxwell-Klein-Gordon system is one of well-studied gauge-field-theoretic wave equations. In $(1+4)$ dimensional setting, the global dynamics of solutions to the system are shown by Oh and Tataru \cite{ohtataru,ohtataru1,ohtataru2}. However, global solutions to the system in $(1+3)$ dimensions is still open. The main drawback of the system is the strong singularity in the quadratic nonlinearity $|\nabla|^{-1}Q_{ij}(\phi,\overline\phi)$. Our first main result provides a partial answer on the question of the scattering property of solutions to the Maxwell-Klein-Gordon system for the scale-invariant Sobolev regularity.
\subsection{Cubic Dirac equations}
Secondly we would like to investigate long-time behaviour of solutions to cubic Dirac equations with the Hartree-type nonlinearty
\begin{align}\label{main-dirac}
\left\{
\begin{array}{l}
	-i\gamma^\mu\partial_\mu\psi+m\psi = [V_b*(\psi^\dagger\psi)]\gamma^0\psi,\\
	\psi|_{t=0} = \psi_0,
\end{array}
\right.	
\end{align}
where $V_b=V_b(x)$ is the Yukawa-type potential given by
$$
V_b(x) = \frac1{4\pi}\frac{e^{-b|x|}}{|x|},\quad b>0,
$$
and $m>0$ is a positive mass.
Recall that we adapt the summation convention. Here $\psi:\mathbb R^{1+3}\rightarrow\mathbb C^4$ is the Dirac spinor field and $\psi^\dagger$ is the complex conjugate transpose of $\psi$, i.e., $\psi^\dagger=(\psi^*)^T$. The Dirac gamma matrices $\gamma^\mu$ are the $4\times4$ complex matrices given by 
\begin{align*}
\gamma^0 = \begin{bmatrix} 
 	I_{2\times2} & \mathbf0 \\ \mathbf0 & -I_{2\times2}
 \end{bmatrix}, \ \gamma^j = \begin{bmatrix} 
 	\mathbf 0 & \sigma^j \\ -\sigma^j & \mathbf0 
 \end{bmatrix}	,
\end{align*}
with the Pauli matrices $\sigma^j$, $j=1,2,3$ given by
\begin{align*}
\sigma^1 = \begin{bmatrix} 
 	0 & 1 \\ 1 & 0 
 \end{bmatrix}, \ \sigma^2 = \begin{bmatrix} 
 	0 & -i \\ i & 0 
 \end{bmatrix}, \  \sigma^3=\begin{bmatrix} 
 	1 & 0 \\ 0 & -1
 \end{bmatrix}.	
\end{align*}
As Theorem \ref{gwp-wave} we prove the global well-posedness and scattering for the scaling critical Sobolev data.
\begin{thm}\label{gwp-dirac}
Let $\sigma=1$. Suppose that the initial data $\psi_0\in L^2_\sigma$ satisfies $\|\psi_0\|_{L^2_\sigma}\ll1$. The Cauchy problem for the equation \eqref{main-dirac} is globally well-posed and scatters to free solutions as $t\rightarrow\pm\infty$.
\end{thm}
\subsubsection{Application to nonlinear Dirac equations}
 Now we shall discuss an application of Theorem \ref{gwp-dirac}. First of all it is instructive to introduce the general form of the Dirac-Klein-Gordon system. Indeed, cubic Dirac equations of the form \eqref{main-dirac} can be obtained by uncoupling the Dirac-Klein-Gordon system 
\begin{align}\label{dkg}
\left\{
\begin{array}{l}
	(-i\gamma^\mu\partial_\mu+M)\psi = g \phi\Gamma\psi, \\
	(\Box+m^2) \phi = -g \psi^\dagger\gamma^0\Gamma\psi.
\end{array}
\right.	
\end{align}
Here $g$ is a coupling constant and we put $g=1$ for simplicity. The $4\times4$ matrix $\Gamma$ can be chosen properly by researchers, for example, $\Gamma = I_{4\times4}, \gamma^0, -\gamma^0\gamma^1\gamma^2\gamma^3$ \cite{bjor}. From \eqref{dkg} one can obtain cubic Dirac equations of the form 
\begin{align}\label{gen-dirac}
(-i\gamma^\mu\partial_\mu+M)\psi = V_b*(\psi^\dagger\gamma^0\Gamma\psi)\Gamma\psi.
\end{align}
We refer the readers to \cite{tes,tes1,cyang} for more detailed derivation from the system \eqref{dkg} to \eqref{gen-dirac}. Recently the nonlinear Dirac systems \eqref{dkg} and \eqref{gen-dirac} with $\Gamma=I_{4\times4}$ have been extensively studied. See \cite{danfos,behe,candyherr,candyherr1,cholee,chohlee,choozlee,wang} and reference therein. For the case $\Gamma=\gamma^\mu$ and the Klein-Gordon field $\phi$ replaced by the vector potential $A_\mu$ with $m=0$, the system \eqref{dkg} becomes the Maxwell-Dirac system \cite{dasfos1,gaoh}. In the case $\Gamma=I_{4\times4}$, it is crucial to exploit the null structure in the bilinear form $\psi^\dagger\gamma^0\psi$ to attain low regularity well-posedness. If $\Gamma=\gamma^0$, however, one cannot enjoy such an advantage and in consequence it is not easy to obtain global well-posedness for a low regularity data. Our second main result says that one can establish scattering property even when it is not possible to take an advantage of null structures.

We would like to mention the Cauchy problems for the boson star equation (or the semi-relativistic equation) with the Hartree-type nonlinearity on $\mathbb R^{1+3}$:
\begin{align}\label{boson-star}
\left\{
\begin{array}{l}
	-i\partial_tu + \sqrt{m^2-\Delta}u = (V_b*|u|^2)u, \\
	u|_{t=0} = u_0
\end{array}\right.	
\end{align}
We refer to \cite{chooz,herrlenz,herrtes} for this well-studied equation. After the use of the Dirac projection operators (see Section \ref{sec:dirac-op}) our Dirac equations \eqref{main-dirac} is of the form \eqref{boson-star}. Thus as a direct application of Theorem \ref{gwp-dirac}, we have 
\begin{cor}\label{dirac-appli}
Let $\sigma=1$. Suppose that the initial data $u_0\in L^2_\sigma$ satisfies $\|u_0\|_{L^2_\sigma}\ll1$. The Cauchy problems for the equation \eqref{boson-star} is globally well-posed and scatters to free solutions as $t\rightarrow\pm\infty$.	
\end{cor}
\noindent By Corollary \ref{dirac-appli} we improve the previous results on the Cauchy problems for \eqref{boson-star} and attain the scaling critical regularity. 

The rest of this paper is organised as follows. In the next section, we give some preliminaries which include half-wave decompositions, Dirac operators, multipliers, definition and basic properties on $U^p-V^p$ spaces and auxiliary estimates. Section 3 and Section 4 are devoted to the proof of our main results, Theorem \ref{gwp-wave} and Theorem \ref{gwp-dirac}, respectively. 
\subsection*{Notations}
\begin{enumerate}
\item
As usual different positive constants, which are independent of dyadic numbers $\mu,\lambda$, and $h$ are denoted by the same letter $C$, if not specified. The inequalities $A \lesssim B$ and $A \gtrsim B$ means that $A \le CB$ and
$A \ge C^{-1}B$, respectively for some $C>0$. By the notation $A \approx B$ we mean that $A \lesssim B$ and $A \gtrsim B$, i.e., $\frac1CB \le A\le CB $ for some absolute constant $C$. We also use the notation $A\ll B$ if $A\le \frac1CB$ for some large constant $C$. Thus for quantities $A$ and $B$, we can consider three cases: $A\approx B$, $A\ll B$ and $A\gg B$. In fact, $A\lesssim B$ means that $A\approx B$ or $A\ll B$.

The spatial and space-time Fourier transform are defined by
$$
\widehat{f}(\xi) = \int_{\mathbb R^3} e^{-ix\cdot\xi}f(x)\,dx, \quad \widetilde{u}(\tau,\xi) = \int_{\mathbb R^{1+3}}e^{-i(t\tau+x\cdot\xi)}u(t,x)\,dtdx.
$$
We also write $\mathcal F_x(f)=\widehat{f}$ and $\mathcal F_{t, x}(u)=\widetilde{u}$. We denote the backward and forward wave propagation of a function $f$ on $\mathbb R^3$ by
$$
e^{-\theta it |\nabla|}f = \frac1{(2\pi)^3}\int_{\mathbb R^3}e^{ix\cdot\xi}e^{-\theta it|\xi|}\widehat{f}(\xi)\,d\xi,
$$
where $\theta\in\{+,-\}$.
\item
We fix a smooth function $\rho\in C^\infty_0(\mathbb R)$ such that $\rho$ is supported in the set $\{ \frac12<r<2\}$ and we let
$$
\sum_{\lambda\in2^{\mathbb Z}}\rho\left(\frac r\lambda\right) =1,
$$
and write $\rho_1=\sum_{\lambda\le1}\rho(\frac r\lambda)$ with $\rho_1(0)=1$.  Now we define the standard Littlewood-Paley multipliers for $\lambda\in 2^{\mathbb N}$ and $\lambda>1$:
$$
P_\lambda = \rho\left(\frac{|-i\nabla|}{\lambda}\right),\quad P_1=\rho_1(|-i\nabla|).
$$
\end{enumerate}

\section{Preliminaries}

\subsection{Half-wave decomposition of the d'Alembertian}
We formulate nonlinear wave equations $\Box u =F$ as a first-order system, which clarifies the dispersive properties of a nonlinear wave. (See also \cite{huhoh}.) We first write
$$
\frac{\partial}{\partial t}\begin{bmatrix} u \\ \partial_tu  \end{bmatrix} = \begin{bmatrix} 0 & 1 \\ \Delta & 0 \end{bmatrix} \begin{bmatrix}u \\ \partial_tu \end{bmatrix} + \begin{bmatrix} 0 \\ F \end{bmatrix}.
$$
We make use of the transform
$$
(u,\partial_tu) \rightarrow (u_+,u_-),\ (0,F) \rightarrow(F_+,F_-),
$$
where
$$
u_\pm  = \frac12\left( u\mp\frac{1}{i|\nabla|}\partial_tu \right),\ F_\pm = \mp\frac{1}{2i|\nabla|}F,
$$
with $|\nabla|=\sqrt{-\Delta}$, which yields the following diagonal system
$$
\frac{\partial}{\partial t} \begin{bmatrix}  u_+ \\ u_- \end{bmatrix} = \begin{bmatrix}  -i|\nabla| & 0 \\ 0 & +i|\nabla| \end{bmatrix} \begin{bmatrix} u_+ \\ u_- \end{bmatrix} + \begin{bmatrix} F_+ \\ F_- \end{bmatrix}. 
$$
This is equivalent to the following half-wave equations
\begin{align}\label{wave-decom}
(-i\partial_t+\theta|\nabla|)u_\theta = \theta\frac{1}{2|\nabla|}F,
\end{align}
where $\theta\in\{+,-\}$. Thus we conclude that the initial value problems for the equation \eqref{main-wave} is reduced to the following first-order system of nonlinear wave equations
\begin{align}
\left\{
\begin{array}{l}
	(-i\partial_t+\theta|\nabla|)u_\theta = \theta|\nabla|^{-2}Q(\overline u,u), \\
	u_{\theta}|_{t=0} = u_{0,\theta}. 
\end{array}
\right.	
\end{align}

\subsection{Dirac projection operators}\label{sec:dirac-op}
Recall the Dirac equations \eqref{main-dirac}
$$
-i\gamma^\mu\partial_\mu\psi + m\psi = [V*(\psi^\dagger\psi)]\gamma^0\psi.
$$ We would like to decompose the Dirac equations and obtain a similar form of a first-order system of half-wave equations as we have done in the previous section. To do this, we first introduce the projections for $\theta\in\{+,-\}$
\begin{align}
\Pi_\theta(\xi) = \frac12\left(I_{4\times4}+\theta\frac{\xi_j\gamma^0\gamma^j+m\gamma^0}{\langle\xi\rangle_m} \right),	
\end{align}
where we used the summation convention and the gamma matrices $\gamma^\mu\in\mathbb C^{4\times4}$, $\mu=0,1,2,3$ are given by 
\begin{align*}
\gamma^0 = \begin{bmatrix} 
 	I_{2\times2} & \mathbf0 \\ \mathbf0 & -I_{2\times2}
 \end{bmatrix}, \ \gamma^j = \begin{bmatrix} 
 	\mathbf 0 & \sigma^j \\ -\sigma^j & \mathbf0 
 \end{bmatrix}	,
\end{align*}
with the Pauli matrices $\sigma^j\in\mathbb C^{2\times2}$, $j=1,2,3$ given by
\begin{align*}
\sigma^1 = \begin{bmatrix} 
 	0 & 1 \\ 1 & 0 
 \end{bmatrix}, \ \sigma^2 = \begin{bmatrix} 
 	0 & -i \\ i & 0 
 \end{bmatrix}, \  \sigma^3=\begin{bmatrix} 
 	1 & 0 \\ 0 & -1
 \end{bmatrix}.	
\end{align*}
Now we define the Fourier multiplier by the identity $\mathcal F_x[\Pi_\theta f](\xi) = \Pi_\theta(\xi)\widehat{f}(\xi)$. By an easy computation one easily see the identity $\Pi_\theta\Pi_\theta=\Pi_\theta$ and $\Pi_\theta\Pi_{-\theta}=0$. We also have $\psi=\Pi_+\psi+\Pi_-\psi$. Then we see that $(-i\gamma^\mu\partial_\mu+m)\Pi_\theta\psi=\gamma^0(-i\partial_t+\theta\langle\nabla\rangle_m)\psi$ and
%\begin{align*}
%(-i\partial_t+\theta\langle\nabla\rangle_m)\Pi_\theta\psi = \Pi_\theta[V*(\psi^\dagger\psi)\psi]. 	
%\end{align*}
hence we conclude that the initial value problems for the equations \eqref{main-dirac} is reduced to the following first-order system of nonlinear Klein-Gordon equations
\begin{align}\label{dirac-decom}
\left\{
\begin{array}{l}
	(-i\partial_t+\theta\langle\nabla\rangle_m)\psi_\theta = \Pi_\theta[V_b*(\psi^\dagger\psi)\psi], \\
	\psi_\theta|_{t=0}=\psi_{0,\theta},
\end{array}
\right.	
\end{align}
 where $\psi_\theta = \Pi_\theta\psi$.
\subsection{Multipliers}\label{multi}
We define $\mathcal Q_\mu$ to be a finitely overlapping collection of cubes of diameter $\frac{\mu}{1000}$ covering $\mathbb R^3$, and let $\{ \rho_{\mathsf q}\}_{\mathsf q\in\mathcal Q_\mu}$ be a corresponding subordinate partition of unity.
For $\mathsf q\in\mathcal Q_\mu$, $d\in 2^{\mathbb Z}$ let
$$
P_{\mathsf q} = \rho_{\mathsf q}(-i\nabla),\quad C^{\theta}_d = \rho\left(\frac{|-i\partial_t + \theta|\nabla||}{d}\right).
$$

We define $C^\theta_{\le d}=\sum_{\delta\le d}C^\theta_\delta$. %and $C^\theta_{\ge d}$ is defined in the similar way. %For simplicity we also write $C^+_h = C_h$.
Given $0 < \alpha \lesssim1$, we define $\mathcal C_{\alpha}$ to be a collection of finitely overlapping caps of radius ${\alpha}$ on the sphere $\mathbb S^2$. If $\kappa\in\mathcal C_{\alpha}$, we let $\omega_\kappa$ be the centre of the cap $\kappa$. Then we define $\{\rho_\kappa\}_{\kappa\in\mathcal C_{\alpha}}$ to be a smooth partition of unity subordinate to the conic sectors $\{ \xi\neq0 , \frac{\xi}{|\xi|}\in\kappa \}$ and denote the angular localisation Fourier multipliers by
$
R_\kappa = \rho_\kappa(-i\nabla).
$

\subsection{Analysis on the sphere}\label{an-sph}
We recall some basic facts from harmonic analysis on the unit sphere. We refer the readers to \cite{candyherr, ster} for the most of ingredients in this section. 
%We also refer the readers to \cite{steinweiss} for more systematic introduction to the spherical harmonics.
We let $Y_{\ell}$ be the set of homogeneous harmonic polynomial of degree $\ell$ on $\mathbb R^3$. Then define $\{ y_{\ell,n} \}_{n=0}^{2\ell}$ a set of orthonormal basis for $Y_{\ell}$, with respect to the inner product:
\begin{align}
\langle y_{\ell,n},y_{\ell',n'}\rangle_{L^2_\omega(\mathbb S^2)} = \int_{\mathbb S^2}{y_{\ell,n}(\omega)} \overline{y_{\ell',n'}(\omega)}\,d\omega.
\end{align}
Given $f\in L^2_x(\mathbb R^3)$, we have the orthogonal decomposition as follow:
\begin{align}
f(x) = \sum_{\ell}\sum_{n=0}^{2\ell}\langle f(|x|\omega),y_{\ell,n}(\omega)\rangle_{L^2_\omega(\mathbb S^2)}y_{\ell,n}\big(\frac{x}{|x|}\big).
\end{align}
For a dyadic number $N>1$, we define the spherical Littlewood-Paley decompositions by
\begin{align}\begin{aligned}\label{hn}
H_N(f)(x) & = 	\sum_{\ell}\sum_{n=0}^{2\ell}\rho\left(\frac\ell N\right)\langle f(|x|\omega),y_{\ell,n}(\omega)\rangle_{L^2_\omega(\mathbb S^2)}y_{\ell,n}\big(\frac{x}{|x|}\big), \\
H_1(f)(x) & = \sum_{\ell}\sum_{n=0}^{2\ell}\rho_{\le1}(\ell)\langle f(|x|\omega),y_{\ell,n}(\omega)\rangle_{L^2_\omega(\mathbb S^2)}y_{\ell,n}\big(\frac{x}{|x|}\big).
\end{aligned}\end{align}
Since $-\Delta_{\mathbb S^2}y_{\ell, n} = \ell(\ell+1)y_{\ell, n}$, by orthogonality one can readily get
$$\|\langle\Omega\rangle^\sigma f\|_{L^2_\omega({\mathbb S^2})} \approx \left\|\sum_{N\in2^{\mathbb N}\cup\{0\}}N^\sigma H_Nf\right\|_{L^2_\omega({\mathbb S^2})}.$$

\begin{lem}[Lemma 7.1 of \cite{candyherr}]\label{sph-ortho}
Let $N\ge1$. Then $H_N$ is uniformly bounded on $L^p(\mathbb R^3)$ in $N$, and $H_N$ commutes with all radial Fourier multipliers. Moreover, if $N'\ge1$, then either $N\approx N'$ or
$$
H_N\Pi_\theta H_{N'}=0.
$$	
%where
%$
%\Pi_{\theta} :=\frac{1}{2}\left(\mathbb I + \theta \Lambda^{-1}\Big[\alpha^x \cdot (-i\nabla) + \beta\Big]\right),
%$
%with $\theta\in\{+,-\}$.
\end{lem}
As an application of Lemma \ref{sph-ortho} one can say that the spherical harmonic projections $H_N$ commutes with the Littlewood-Paley projections such as $P_\lambda$ and $C^\theta_d$. Furtheremore the orthogonality of the spherical harmonics still holds when one deals with the Dirac projections. 
\subsection{Adapted function spaces}\label{ftn-sp}
We discuss the basic properties of function spaces of $U^p$ and $V^p$ type. We refer the readers to \cite{haheko,kochtavi} for more details. Let $\mathcal I$ be the set of finite partitions $-\infty=t_0<t_1<\cdots<t_K=\infty$ and let $1\le p<\infty$.
\begin{defn}
A function $a:\mathbb R\rightarrow L^2_x$ is called a $U^p$-atom if there exists a decomposition 
$$
a=\sum_{j=1}^K\chi_{[t_k-1,t_k)}(t)f_{j-1}
$$
 with
$$
\{f_j\}_{j=0}^{K-1}\subset L^2_x,\ \sum_{j=0}^{K-1}\|f_j\|_{L^2_x}^p=1,\ f_0=0.
$$
Furthermore, we define the atomic Banach space 
$$
U^p := \left\{ u=\sum_{j=1}^\infty \lambda_ja_j : a_j \, U^p\textrm{-atom},\ \lambda_j\in\mathbb C \textrm{ such that } \sum_{j=1}^\infty|\lambda_j|<\infty  \right\}
$$
with the induced norm
$$
\|u\|_{U^p} := \inf\left\{ \sum_{j=1}^\infty|\lambda_j| : u=\sum_{j=1}^\infty \lambda_ja_j,\,\lambda_j\in\mathbb C,\, a_j \, U^p\textrm{-atom}   \right\}.
$$
\end{defn}
We list some basic properties of $U^p$ spaces. 
\begin{prop}[Proposition 2.2 of \cite{haheko}]
Let $1\le p<q<\infty$.
\begin{enumerate}
\item $U^p$ is a Banach space.
\item The embeddings $U^p\subset U^q\subset L^\infty(\mathbb R;L^2_x)$ are continuous.
\item For $u\in U^p$, $u$ is right-continuous.
\end{enumerate}
\end{prop}
\noindent We also define the space $U^p_\theta$ to be the set of all $u\in\mathbb R\rightarrow L^2_x$ such that $e^{-\theta it|\nabla|}u\in U^p$ with the obvious norm
$
\|u\|_{U^p_\theta} := \|e^{-\theta it|\nabla|}u\|_{U^p}.
$
%Let $\mathcal Z =\left\{ \{t_k\}_{k=0}^K : t_k\in\mathbb R, t_k<t_{k+1} \right\}$ be the set of increasing sequences of real numbers.
We define the $2$-variation of $v$ to be
$$
|v|_{V^2} = \sup_{ \{t_k\}_{k=0}^K\in\mathcal I } \left( \sum_{k=0}^K\|v(t_k)-v(t_{k-1})\|_{L^2_x}^2 \right)^\frac12
$$
Then the Banach space $V^2$ can be defined to be all right continuous functions $v:\mathbb R\rightarrow L^2_x$ such that the quantity
$$
\|v\|_{V^2} = \|v\|_{L^\infty_tL^2_x} + |v|_{V^2}
$$
is finite. Set $\|u\|_{V^2_\theta}=\|e^{-\theta it|\nabla|}u\|_{V^2}$. We recall basic properties of $V^2_\theta$ space from \cite{candyherr, candyherr1, haheko}. In particular,  we use the following lemma to prove the scattering result.
\begin{lem}[Lemma 7.4 of \cite{candyherr}]\label{v-scatter}
	Let $u\in V^2_\theta$. Then there exists $f\in L^2_x$ such that $\|u(t)-e^{-\theta it|\nabla|}f\|_{L^2_x}\rightarrow0$ as $t\rightarrow\pm\infty$.
\end{lem}
The following lemma is on a simple bound in the high-modulation region.
\begin{lem}[Corollary 2.18 of \cite{haheko}]
Let $2\le q\le\infty$. For $d\in2^{\mathbb Z}$ and $\theta \in \{+, -\}$, we have
\begin{align}\label{bdd-high-mod}
\begin{aligned}
\|C^{\theta}_d u\|_{L^q_tL^2_x} \lesssim d^{-\frac1q}\|u\|_{V^2_\theta},
%\|C^{\theta}_{\ge h}u\|_{L^q_tL^2_x} \lesssim h^{-\frac1q}\|u\|_{V^2_\theta}.
\end{aligned}
\end{align}
\end{lem}
\begin{lem}[Lemma 2.2 of \cite{cholee}]\label{energy-ineq}
Let $u\in U^p$ be absolutely continuous with $1<p<\infty$. Then 
\begin{align*}
\|u\|_{U^p} = \sup\left\{ \left| \int \langle u'(t),v(t)\rangle_{L^2_x}\,dt \right| : v\in C^\infty_0,\ \|v\|_{V^{p'}}=1 \right\}.	
\end{align*}

%$$
%\sup_{\|P_\lambda H_Nw\|_{V^2_\theta}\le1}\left| \int_{\mathbb R}\langle P_\lambda H_N w(t),F(t)\rangle_{L^2_x}\,dt \right|<\infty.
%$$	
%If $u\in C(\mathbb R; L^2_x)$ satisfies $-i\partial_t u+\theta|\nabla| u=F$, then $P_\lambda H_N u\in V^2_\theta$ and we have the bound
%\begin{align}
%\|P_\lambda H_N u\|_{V^2_\theta} \lesssim \|P_\lambda H_N u(0)\|_{L^2_x}+\sup_{\|P_\lambda H_Nw\|_{V^2_\theta}\le1}\left| \int_{\mathbb R}\langle P_\lambda H_N w(t),F(t)\rangle_{L^2_x}\,dt \right|.	
%\end{align}
\end{lem}
We define the Banach space associated with the homogeneous Sobolev space to be the set
$$
\dot F^{s,\sigma}_\theta = \left\{ u\in C(\mathbb R;\langle\Omega\rangle^{-\sigma}\dot H^s): \|u\|_{\dot F^{s,\sigma}_\theta}<\infty \right\},
$$
where the norm is defined by
$$
\|u\|_{\dot F^{s,\sigma}_\theta} = \bigg( \sum_{\lambda\in2^{\mathbb Z}}\sum_{N\ge1}\lambda^{2s}N^{2\sigma}\|P_\lambda H_Nu\|_{U^2_{\theta}}^2 \bigg)^\frac12.
$$
Similarly we define the Banach space $F^{s,\sigma}_\theta$ associated to the inhomogenous Sobolev space in the obvious way. 
\begin{rem}
So far we have defined the adapted function spaces for the wave operator. However, with a slight modification we can also define the adapted function spaces for the Klein-Gordon-type operator and all the above lemma also holds for the Klein-Gordon operator. In consequence, for brevity we allow abuse of notation and simply use the notation $U^p_\theta$ and $V^p_\theta$ for both wave and Klein-Gordon operators. We also refer the readers to \cite{candyherr}.  
\end{rem}

\subsection{Auxiliary estimates}
We begin with very basic Sobolev estimates which is also known as the Bernstein inequality. 
\begin{lem}
Let $0<\alpha\lesssim1$ and $\kappa\in\mathcal C_\alpha$. Let $\lambda>0$ be a dyadic number. For any test function $f$ on $\mathbb R^3$ we have
\begin{align}\label{bernstein}
\|R_\kappa P_\lambda f\|_{L^\infty_x} \lesssim (\lambda^3\alpha^2)^\frac1p \|f\|_{L^p_x}.  
\end{align}
\end{lem}
To study the dispersive property of solutions it is of great importance to exploit so-called the Strichartz estimates \cite{keeltao,stri}. In this paper we use an improved Strichartz estimate which is obtained by spending an extra regularity with respect to the angular variables. (See also \cite{choslee,ster}.)
\begin{prop}
For $\frac{1}{10}\ge\eta>0$, let $q_\eta=\frac{4}{1-\eta}$. We have the improved Strichartz estimates by imposing angular regularity as follow:
\begin{align}\label{stri-ang}
\|e^{\theta it|\nabla|}P_\lambda H_N f\|_{L^2_tL^{q_\eta}_x} \lesssim \lambda^{1-\frac3{q_\eta}}N^{\frac12+\eta}\|	P_\lambda H_Nf\|_{L^2_x}.
\end{align}
%Similarly, we have
%\begin{align}
%\|e^{it\Lambda}P_\lambda H_N f\|_{L^q_tL^4_x} \lesssim \lambda^{\frac34-\frac1q}N\|P_\lambda H_Nf\|_{L^2_x},	
%\end{align}
%provided that $2<q\le\infty$.
\end{prop}
\noindent %Throughout this paper we shall simply denote $q_\eta$ by $q$. 
Note that the estimates hold when we replace the propogator $e^{-it|\nabla|}$ with $e^{-it\langle\nabla\rangle}$. The space-time estimates \eqref{stri-ang} say that one can obtain an improved bound when dealing with multilinear estimates. However, the singularity $|\nabla|^{-2}$ in the nonlinearity in \eqref{main-wave} is too strong and we cannot obtain the desired well-posedness for the critical Sobolev data by simply using the estimate \eqref{stri-ang}. This is why the low-output frequency interaction becomes the most serious case. To overcome this problem, we apply the almost orthogonal decompositions of conic sectors. The question is whether one can obtain the better estimates by exploiting the localisation into the conic sectors. The following lemma which is also known as {\it angular concentration estimates} answers this question.
\begin{lem}[Lemma 8.5 of \cite{candyherr}]\label{ang-con}
Let $2\le p<\infty$, and $0\le s<\frac2p$. If $\lambda,N\ge1$, ${\alpha}\gtrsim\lambda^{-1}$, and $\kappa\in\mathcal C_{\alpha}$, then we have
$$
\|R_\kappa P_\lambda H_N f\|_{L^p_x(\mathbb R^3)} \lesssim ({\alpha} N)^s \|P_\lambda H_N f\|_{L^p_x(\mathbb R^3)}.
$$	
\end{lem}
\noindent We refer the readers to \cite{sterbenz2} for the proof. Note that in the Bernstein inequality, it is no harm to put $f\rightarrow R_{\kappa'}P_{2\lambda}f$ with $\kappa'\in\mathcal C_{2\alpha}$. With cube localisation $P_\mathtt q$ of size $\mu\le\lambda$, we use in order the Bernstein inequality, angular concentration estimates, and then the improved Strichartz estimates. Here one should note that $\frac12-2\eta<\frac2{q_\eta}=\frac{1-\eta}{2}$. Then we see that
\begin{align*}
\|P_\mathtt q R_\kappa P_\lambda H_N u\|_{L^2_tL^\infty_x} & \lesssim (\mu^3\alpha^2)^\frac1q \sup_{\kappa\in\mathcal C_\alpha}\|R_\kappa P_\lambda H_N u\|_{L^2_tL^q_x} \\
& \lesssim 	(\mu^3\alpha^2)^\frac1q(\alpha N)^{\frac12-2\eta}\|P_\lambda H_Nu\|_{L^2_tL^q_x} \\
& \lesssim \mu^\frac3q\alpha^\frac2q\alpha^{\frac12-2\eta}N^{\frac12-2\eta}\lambda^{1-\frac3q}N^{\frac12+\eta}\|P_\lambda H_N u\|_{U^2_\theta}.
\end{align*}
The above argument will be often used in the proof of Theorem \ref{gwp-wave} and Theorem \ref{gwp-dirac}.
\section{Bilinear estimates: Proof of Theorem \ref{gwp-wave}}
Now we arrive at the proof of Theorem \ref{gwp-wave}. First we define the Duhamel integral
\begin{align*}
\mathfrak I^\theta[F] = \int_0^t e^{-\theta i(t-t')|\nabla|}F(t')\,dt'.	
\end{align*}
Then $\mathfrak I^\theta[F]$ solves the equation 
$$
(-i\partial_t+\theta|\nabla|)\mathfrak I^\theta[F] = F,
$$
with vanishing data at $t=0$. 
To prove Theorem \ref{gwp-wave} it is enough to show the following bilinear estimates
\begin{align}
\left\|\mathfrak I^\theta[|\nabla|^{-2}Q(\overline u,v)]\right\|_{\dot F^{\frac12,1}_\theta} & \lesssim \|u\|_{\dot F^{\frac12,1}_{\theta_1}}\|v\|_{\dot F^{\frac12,1}_{\theta_2}}. \label{main-bi}
%\left\|\mathfrak I^\theta[\Lambda^{-1}Q(\Lambda^{-1}u,v)]\right\|_{F^{\frac12,1}_\theta} & \lesssim \|u\|_{F^{\frac12,1}_{\theta_1}}\|v\|_{F^{\frac12,1}_{\theta_2}}.
\end{align}
Then an application of the standard contraction argument gives the desired global solutions to the equation \eqref{main-wave} when we have the smallness assumptions on the initial datum:
$
\|(u_0,u_1)\|_{\dot H^\frac12_\sigma\times\dot H^{-\frac12}_\sigma} \ll1.
$
Moreover, the continuous embedding $U^2\subset V^2$ and Lemma \ref{v-scatter} imply the scattering in $U^2_\theta$ space.
By the duality in $U^2-V^2$ Lemma \ref{energy-ineq} we obtain the trilinear expression as follows
\begin{align*}
	\left\|\mathfrak I^\theta[|\nabla|^{-2}Q(\overline u,v)]\right\|^2_{\dot F^{\frac12,1}_\theta} & \lesssim \sum_{\mu\in2^\mathbb Z}\sum_{N\ge1}(\mu^\frac12 N)^2 \|P_\mu H_N \mathfrak I^\theta[|\nabla|^{-2}Q(\overline u,v)]\|_{U^2_\theta}^2 \\
	& \lesssim \sum_{\mu\in2^\mathbb Z}\sum_{N\ge1}(\mu^\frac12N)^2\sup_{\|P_\mu H_Nw\|_{V^2_\theta}\le1}\left| \int_{\mathbb R^{3+1}}P_\mu H_Nw\, |\nabla|^{-2}Q(\overline u,v) \,dxdt \right|^2.
\end{align*}
Thus our main bilinear estimates can be obtained provided that the following frequency-localised trilinear estimates holds:
\begin{lem}
Let $0<\eta\ll1$ be a small positive number. For some $\frac18<\delta\le\frac14$, we have
\begin{align}\label{main-bi-loc}
\begin{aligned}
&\left|\int_{\mathbb R^{1+3}} w_{\mu,N}\,  |\nabla|^{-2}Q(\overline{u_{\lambda_1,N_1}},v_{\lambda_2,N_2}) \,dxdt\right| \\
 & \qquad\qquad \lesssim (\min\{\lambda_1,\lambda_2\})^\frac12 \left( \frac{\min\{\mu,\lambda_1,\lambda_2\}}{\max\{\mu,\lambda_1,\lambda_2\}} \right)^\delta (\min\{N_1,N_2\})^{1-\eta}	 \|w_{\mu,N}\|_{V^2_\theta}\|u_{\lambda_1,N_1}\|_{U^2_{\theta_1}}\|v_{\lambda_2,N_2}\|_{U^2_{\theta_2}},
 \end{aligned}
\end{align}
\end{lem}
\noindent where we put $w_{\mu,N}=P_\mu H_Nw, u_{\lambda_1,N_1}=P_{\lambda_1}H_{N_1}u,$ and $v_{\lambda_2,N_2}=P_{\lambda_2}H_{N_2}v$ for brevity. To obtain \eqref{main-bi-loc} we shall consider all possible frequency interactions.
In view of the standard Littlewood-Paley trichotomy one can easily see that the integral in \eqref{main-bi-loc} vanishes unless the following interactions hold:
\begin{align}
\min\{\mu,\lambda_1,\lambda_2\}	& \lesssim \textrm{med}\{\mu,\lambda_1,\lambda_2\} \approx \max\{\mu,\lambda_1,\lambda_2\}, \\
\min\{N,N_1,N_2\} & \lesssim \textrm{med}\{N,N_1,N_2\} \approx \max \{N,N_1,N_2\}.
\end{align}
We first decompose the integrand in \eqref{main-bi-loc} with respect to the modulation as follows
\begin{align*}
	&\int_{\mathbb R^{1+3}} P_\mu H_Nw\,  |\nabla|^{-2}Q(\overline{u_{\lambda_1,N_1}},v_{\lambda_2,N_2}) \,dxdt \\
	 & = \sum_{d\in2^{\mathbb Z}}\int_{\mathbb R^{3+1}}C^\theta_d w_{\mu,N}|\nabla|^{-2}Q(C^{\theta_1}_{\ll d}\overline{u_{\lambda_1,N_1}},C^{\theta_2}_{\ll d}v_{\lambda_2,N_2})\,dtdx \\
	 & \qquad + \sum_{d\in2^{\mathbb Z}}\int_{\mathbb R^{3+1}}C^\theta_{\le d}w_{\mu,N}|\nabla|^{-2}Q(C^{\theta_1}_d\overline{u_{\lambda_1,N_1}},C^{\theta_2}_{\le d}v_{\lambda_2,N_2})\,dtdx \\
	 & \qquad + \sum_{d\in2^{\mathbb Z}}\int_{\mathbb R^{3+1}}C^\theta_{\le d}w_{\mu,N}|\nabla|^{-2}Q(C^{\theta_1}_{\le d}\overline{u_{\lambda_1,N_1}},C^{\theta_2}_{ d}v_{\lambda_2,N_2})\,dtdx \\
	 & := \sum_{d\in2^{\mathbb Z}} \mathcal I_0+\mathcal I_1+\mathcal I_2.
\end{align*}
\subsection{Low modulation}\label{sec:low-mod}
Now we consider the low-modulation regime
$
d \lesssim \min\{\mu,\lambda_1,\lambda_2\}.
$
In this regime we will pay special attention  to the low-output interaction, i.e., $\mu\ll\lambda_1\approx\lambda_2$. The main problem is that even when we can take an advantage of the presence of null forms very favourably in the low-output interactions, the Fourier multiplier $|\nabla|^{-2}$ gives rise to the serious singularity and the cancellation property given by null structure is not sufficient to cover all such a {\it bad interaction}. To overcome this problem we adapt fully angular momentum operator and exploit angular concentration phenomena via bilinear decompositions by conic sectors. The key point is that when one input-frequency is localised in a conic sector of a small angle, the other input-frequency should be also localised in another conic sector of a compatible size. On the other hand, in the high-output interaction the null structure no longer plays any crucial role compared to the low-output case since we only have bilinear decompositions by a rather wide-angle. This is not problematic however, the Fourier multiplier $|\nabla|^{-2}$ no longer is a serious singularity, instead, it plays a crucial role as strong decay. In consequence, the high-output interaction becomes the easiest case in the proof. 
\subsubsection{Case 1: $\mu\ll \lambda_1\approx\lambda_2$} It is no harm to put $\lambda_1=\lambda_2=\lambda$ in our argument. Put $\alpha = (\frac{d\mu}{\lambda^2})^\frac12$. We first use an almost orthogonal decomposition by smaller cubes and angular sectors and obtain
\begin{align*}
\mathcal I_0 & \lesssim \sum_{\substack{ \mathtt q_1,\mathtt q_2\in \mathcal Q_\mu \\ |\mathtt q_1-\mathtt q_2|\lesssim\mu  }}\sum_{\substack{ \kappa_1,\kappa_2\in\mathcal C_\alpha \\ |\kappa_1-\kappa_2|\lesssim\alpha }}\int_{\mathbb R^{3+1}}C^\theta_d w_{\mu,N}|\nabla|^{-2}Q(P_{\mathtt q_1}R_{\kappa_1}C^{\theta_1}_{\ll d}u_{\lambda,N_1},P_{\mathtt q_2}R_{\kappa_2}C^{\theta_2}_{\ll d}v_{\lambda,N_2})\,dtdx \\
& \lesssim 	\sum_{\substack{ \mathtt q_1,\mathtt q_2\in \mathcal Q_\mu \\ |\mathtt q_1-\mathtt q_2|\lesssim\mu  }}\sum_{\substack{ \kappa_1,\kappa_2\in\mathcal C_\alpha \\ |\kappa_1-\kappa_2|\lesssim\alpha }}\|C^\theta_d w_{\mu,N}\|_{L^2_{t}L^2_x}\||\nabla|^{-2}Q(P_{\mathtt q_1}R_{\kappa_1}C^{\theta_1}_{\ll d}u_{\lambda,N_1},P_{\mathtt q_2}R_{\kappa_2}C^{\theta_2}_{\ll d}v_{\lambda,N_2})\|_{L^2_tL^2_x} \\
& \lesssim d^{-\frac12}\|w_{\mu,N}\|_{V^2_\theta} 	\bigg(\sum_{\substack{ \mathtt q_1,\mathtt q_2\in \mathcal Q_\mu \\ |\mathtt q_1-\mathtt q_2|\lesssim\mu  }}\sum_{\substack{ \kappa_1,\kappa_2\in\mathcal C_\alpha \\ |\kappa_1-\kappa_2|\lesssim\alpha }}\||\nabla|^{-2}Q(P_{\mathtt q_1}R_{\kappa_1}C^{\theta_1}_{\ll d}u_{\lambda,N_1},P_{\mathtt q_2}R_{\kappa_2}C^{\theta_2}_{\ll d}v_{\lambda,N_2})\|_{L^2_tL^2_x}^2\bigg)^\frac12, 
\end{align*}
where we used the simple bound for a high-modulation-regime \eqref{bdd-high-mod} for $w_{\mu,N}$. Now we exploit the null structure in the bilinear form $Q$ and then use the H\"older inequality and Bernstein inequality for $u_{\lambda,N_1}$. In sequel we put $q=\frac{4}{1-\eta}$ for a small $\eta>0$. Then we have
\begin{align*}
\mathcal I_0  & \lesssim d^{-\frac12} \mu^{-2}\alpha \lambda^2 \|w_{\mu,N}\|_{V^2_\theta} \bigg(\sum_{\substack{ \mathtt q_1,\mathtt q_2\in \mathcal Q_\mu \\ |\mathtt q_1-\mathtt q_2|\lesssim\mu  }}\sum_{\substack{ \kappa_1,\kappa_2\in\mathcal C_\alpha \\ |\kappa_1-\kappa_2|\lesssim\alpha }}\|P_{\mathtt q_1}R_{\kappa_1}C^{\theta_1}_{\ll d}u_{\lambda,N_1}\|_{L^2_tL^\infty_x}^2\|P_{\mathtt q_2}R_{\kappa_2}C^{\theta_2}_{\ll d}v_{\lambda,N_2}\|_{L^\infty_tL^2_x}^2\bigg)^\frac12 \\
& \lesssim d^{-\frac12} \mu^{-2}\alpha \lambda^2 \mu^\frac3q\alpha^\frac2q \|w_{\mu,N}\|_{V^2_\theta}\sup_{\kappa_1\in\mathcal C_\alpha}\|R_{\kappa_1}C^{\theta_1}_{\ll d}u_{\lambda,N_1}\|_{L^2_tL^q_x} \\
& \qquad\qquad \times \bigg(\sum_{\substack{ \mathtt q_1,\mathtt q_2\in \mathcal Q_\mu \\ |\mathtt q_1-\mathtt q_2|\lesssim\mu  }}\sum_{\substack{ \kappa_1,\kappa_2\in\mathcal C_\alpha \\ |\kappa_1-\kappa_2|\lesssim\alpha }}\|P_{\mathtt q_2}R_{\kappa_2}C^{\theta_2}_{\ll d}v_{\lambda,N_2}\|_{L^\infty_tL^2_x}^2\bigg)^\frac12.
\end{align*}
The final step is an application of the angular concentration estimates Lemma \ref{ang-con} with $s=\frac12-2\eta<\frac2q$ and then the improve Strichartz estimates \eqref{stri-ang}, which gives
\begin{align*}
\mathcal I_0 & \lesssim d^{-\frac12} \mu^{-2}\alpha \lambda^2 \mu^\frac3q\alpha^\frac2q(\alpha N_1)^{\frac12-2\eta} \|w_{\mu,N}\|_{V^2_\theta} \|C_{\ll d}^{\theta_1}u_{\lambda,N_1}\|_{L^2_tL^q_x}\|C^{\theta_2}_{\ll d}v_{\lambda,N_2}\|_{L^\infty_tL^2_x} \\
& \lesssim d^{-\frac12} \mu^{-2}\alpha \lambda^2 \mu^\frac3q\alpha^\frac2q(\alpha N_1)^{\frac12-2\eta}\lambda^{1-\frac3q} N_1^{\frac12+\eta} \|w_{\mu,N}\|_{V^2_\theta} \|u_{\lambda,N_1}\|_{U^2_{\theta_1}}\|v_{\lambda,N_2}\|_{U^2_{\theta_2}}.
\end{align*}
The summation with respect to $d\lesssim\mu$ yields 
$$
\sum_{d:d\lesssim\mu}\mathcal I_0 \lesssim \mu^{-1+\frac5q-2\eta}\lambda^{\frac32-\frac5q+2\eta}N_1^{1-\eta}\|w_{\mu,N}\|_{V^2_\theta} \|u_{\lambda_1,N_1}\|_{U^2_{\theta_1}}\|v_{\lambda_2,N_2}\|_{U^2_{\theta_2}}.
$$
If $N_1\gg N_2$, then we simply interchange the role of $u_{\lambda,N_1}$ and $v_{\lambda,N_2}$ and then obtain exactly the same bound. We now consider $\mathcal I_1$. As we have done in the previous estimate, we use an almost orthogonal decomposition of cubes and angular sectors to get
\begin{align*}
\mathcal I_1 & \lesssim 	\sum_{\substack{ \mathtt q_1,\mathtt q_2\in \mathcal Q_\mu \\ |\mathtt q_1-\mathtt q_2|\lesssim\mu  }}\sum_{\substack{ \kappa_1,\kappa_2\in\mathcal C_\alpha \\ |\kappa_1-\kappa_2|\lesssim\alpha }}\int_{\mathbb R^{3+1}}C^\theta_{\le d} w_{\mu,N}|\nabla|^{-2}Q(P_{\mathtt q_1}R_{\kappa_1}C^{\theta_1}_{ d}u_{\lambda,N_1},P_{\mathtt q_2}R_{\kappa_2}C^{\theta_2}_{\le d}v_{\lambda,N_2})\,dtdx. 
\end{align*}
The next step is to exploit the null structire and use the H\"older inequality as the previous estimate
\begin{align*}
\mathcal I_1 & \lesssim \mu^{-2}\lambda^2\alpha 	\sum_{\substack{ \mathtt q_1,\mathtt q_2\in \mathcal Q_\mu \\ |\mathtt q_1-\mathtt q_2|\lesssim\mu  }}\sum_{\substack{ \kappa_1,\kappa_2\in\mathcal C_\alpha \\ |\kappa_1-\kappa_2|\lesssim\alpha }}  \int_{\mathbb R^{3+1}} C^\theta_{\le d} w_{\mu,N}P_{\mathtt q_1}R_{\kappa_1}C^{\theta_1}_{ d}u_{\lambda,N_1}P_{\mathtt q_2}R_{\kappa_2}C^{\theta_2}_{\le d}v_{\lambda,N_2}\,dtdx \\
& \lesssim \mu^{-2}\lambda^2\alpha 	\sum_{\substack{ \mathtt q_1,\mathtt q_2\in \mathcal Q_\mu \\ |\mathtt q_1-\mathtt q_2|\lesssim\mu  }}\sum_{\substack{ \kappa_1,\kappa_2\in\mathcal C_\alpha \\ |\kappa_1-\kappa_2|\lesssim\alpha }}  \|  C^\theta_{\le d} w_{\mu,N}\|_{L^\infty_tL^2_x}\| P_{\mathtt q_1}R_{\kappa_1}C^{\theta_1}_{ d}u_{\lambda,N_1}\|_{L^2_tL^2_x}\|P_{\mathtt q_2}R_{\kappa_2}C^{\theta_2}_{\le d}v_{\lambda,N_2}\|_{L^2_tL^\infty_x} \\
& \lesssim \mu^{-2}\lambda^2\alpha \|  C^\theta_{\le d} w_{\mu,N}\|_{L^\infty_tL^2_x}  \bigg(\sum_{\substack{ \mathtt q_1,\mathtt q_2\in \mathcal Q_\mu \\ |\mathtt q_1-\mathtt q_2|\lesssim\mu  }}\sum_{\substack{ \kappa_1,\kappa_2\in\mathcal C_\alpha \\ |\kappa_1-\kappa_2|\lesssim\alpha }} \| P_{\mathtt q_1}R_{\kappa_1}C^{\theta_1}_{ d}u_{\lambda,N_1}\|_{L^2_tL^2_x}^2\|P_{\mathtt q_2}R_{\kappa_2}C^{\theta_2}_{\le d}v_{\lambda,N_2}\|_{L^2_tL^\infty_x}^2\bigg)^\frac12.
\end{align*}
Then we use the Bernstein inequality for $v_{\lambda,N_2}$ and then Lemma \ref{ang-con} and the Strichartz estimates \eqref{stri-ang} 
\begin{align*}
\mathcal I_1 %& \lesssim \mu^{-2}\lambda^2\alpha \mu^\frac3q\alpha^\frac2q \|w_{\mu,N}\|_{V^2_\theta}  \bigg(\sum_{\substack{ \mathtt q_1,\mathtt q_2\in \mathcal Q_\mu \\ |\mathtt q_1+\mathtt q_2|\lesssim\mu  }}\sum_{\substack{ \kappa_1,\kappa_2\in\mathcal C_\alpha \\ |\kappa_1+\kappa_2|\lesssim\alpha }} \| P_{\mathtt q_1}R_{\kappa_1}C^{\theta_1}_{ d}u_{\lambda_1,N_1}\|_{L^2_tL^2_x}^2\|P_{\mathtt q_2}R_{\kappa_2}C^{\theta_2}_{\le d}v_{\lambda_2,N_2}\|_{L^2_tL^q_x}^2\bigg)^\frac12 \\
& \lesssim \mu^{-2}\lambda^2\alpha \mu^\frac3q\alpha^\frac2q \|w_{\mu,N}\|_{V^2_\theta} \sup_{\kappa_2\in\mathcal C_\alpha}\|R_{\kappa_2}C^{\theta_2}_{\le d}v_{\lambda,N_2}\|_{L^2_tL^q_x}  \bigg(\sum_{\substack{ \mathtt q_1,\mathtt q_2\in \mathcal Q_\mu \\ |\mathtt q_1-\mathtt q_2|\lesssim\mu  }}\sum_{\substack{ \kappa_1,\kappa_2\in\mathcal C_\alpha \\ |\kappa_1-\kappa_2|\lesssim\alpha }} \| P_{\mathtt q_1}R_{\kappa_1}C^{\theta_1}_{ d}u_{\lambda,N_1}\|_{L^2_tL^2_x}^2\bigg)^\frac12 \\
& \lesssim \mu^{-2}\lambda^2\alpha \mu^\frac3q\alpha^\frac2q (\alpha N_2)^{\frac12-2\eta} \|w_{\mu,N}\|_{V^2_\theta}  \|C^{\theta_2}_{\le d}v_{\lambda_2,N_2}\|_{L^2_tL^q_x}\|C^{\theta_1}_d u_{\lambda_1,N_1}\|_{L^2_tL^2_x} \\
& \lesssim \mu^{-2}\lambda^2\alpha \mu^\frac3q\alpha^\frac2q (\alpha N_2)^{\frac12-2\eta} \lambda^{1-\frac3q}N_2^{\frac12+\eta}d^{-\frac12}\|w_{\mu,N}\|_{V^2_\theta} \|u_{\lambda_1,N_1}\|_{U^2_{\theta_1}}\|v_{\lambda_2,N_2}\|_{U^2_{\theta_2}},
\end{align*}
where we used the bound \eqref{bdd-high-mod} for $C^\theta_du$. The summation with respect to the modulation $d\lesssim\mu$ gives the desired bound. If $N_1\ll N_2$, we can simply interchange the role of $u_{\lambda,N_1}$ and $v_{\lambda,N_2}$ and follow the above argument. The estimate of $\mathcal I_2$ can be obtained in the identical manner as the estimate of $\mathcal I_1$. We omit the details.
\subsubsection{Case 2: $\lambda_1 \lesssim \mu\approx\lambda_2$} The case $\lambda_2\lesssim\mu\approx\lambda_1$ would readily follow by symmetry and we focus on the case $\lambda_1\ll \lambda_2$. The high-output case is much easier than the low-output case, i.e., $\min\{\mu,\lambda_1,\lambda_2\}=\mu$, since the Fourier multiplier $|\nabla|^{-2}$ in the integrand is not the serious singularity, even further it plays a role as a strong decay. We only treat the estimate of $\mathcal I_1$ with $N_1\gg N_2$ in this paper, since this case is the most serious interaction in the high-output interaction. 
We put $\beta = (\frac d{\lambda_1})^\frac12$ and use the orthogonal decompositions
\begin{align*}
\mathcal I_1 & \lesssim \sum_{\substack{\mathtt q,\mathtt q_2\in\mathcal Q_{\lambda_1} \\ |\mathtt q+\theta_2\mathtt q_2|\lesssim\lambda_1}}\sum_{\substack{\kappa,\kappa_1,\kappa_2\in\mathcal C_\beta \\ |\kappa_1+\kappa_2|,|\kappa+\theta_2\kappa_2|\lesssim\beta}}\int_{\mathbb R^{1+3}}P_\mathtt q R_\kappa C^\theta_{\le d}|\nabla|^{-2}Q(R_{\kappa_1}C^{\theta_1}_du_{\lambda_1,N_1},P_{\mathtt q_2}R_{\kappa_2}C^{\theta_2}_{\le d}v_{\lambda_2,N_2})\,dtdx \\
& \lesssim \mu^{-2}\lambda_1\lambda_2\beta	\sum_{\substack{\mathtt q,\mathtt q_2\in\mathcal Q_{\lambda_1} \\ |\mathtt q+\theta_2\mathtt q_2|\lesssim\lambda_1}}\sum_{\substack{\kappa,\kappa_1,\kappa_2\in\mathcal C_\beta \\ |\kappa_1+\kappa_2|,|\kappa+\theta_2\kappa_2|\lesssim\beta}}\int_{\mathbb R^{1+3}}P_\mathtt q R_\kappa C^\theta_{\le d}R_{\kappa_1}C^{\theta_1}_du_{\lambda_1,N_1}P_{\mathtt q_2}R_{\kappa_2}C^{\theta_2}_{\le d}v_{\lambda_2,N_2}\,dtdx.
\end{align*}
Then we use the H\"older inequality and then the Cauchy-Schwarz inequality in $\kappa_1$ to get
\begin{align*}
\mathcal I_1 & \lesssim 	\mu^{-2}\lambda_1\lambda_2\beta \bigg( \sum_{\kappa_1}\|R_{\kappa_1} C^{\theta_1}_d u_{\lambda_1,N_1}\|_{L^2_tL^2_x}^2 \bigg)^\frac12 \\
& \qquad\times \bigg( \sum_{\kappa_1} \bigg( \sum_{\kappa,\kappa_2}\sum_{\mathtt q,\mathtt q_2}\|P_\mathtt qR_\kappa C^{\theta}_{\le d}w_{\mu,N}\|_{L^\infty_tL^2_x}\|P_{\mathtt q_2}R_{\kappa_2}C^{\theta_2}_{\le d}v_{\lambda_2,N_2}\|_{L^2_tL^\infty_x} \bigg)^2 \bigg)^\frac12.
\end{align*}
We use the Bernstein inequality for $v_{\lambda_2,N_2}$ and obtain
\begin{align*}
\mathcal I_1 & \lesssim 		\mu^{-2}\lambda_1\lambda_2\beta \lambda_1^\frac3q \beta^\frac2q \|C^\theta_du_{\lambda_1,N_1}\|_{L^2_tL^2_x} \sup_{\kappa_2}\|R_{\kappa_2}C^{\theta_2}_{\le d}v_{\lambda_2,N_2}\|_{L^2_tL^q_x} \\
& \qquad\qquad \times\bigg( \sum_{\kappa,\kappa_1,\kappa_2}\sum_{\mathtt q,\mathtt q_2} \|P_\mathtt qR_\kappa C^{\theta}_{\le d}w_{\mu,N}\|_{L^\infty_tL^2_x}^2\bigg)^\frac12.
\end{align*}
The remaining step is to apply the bound for the high-modulation-region \eqref{bdd-high-mod} for $C^\theta_du$ and then Lemma \ref{ang-con} followed by the Strichartz estimate \eqref{stri-ang} for $v_{\lambda_2,N_2}$
\begin{align*}
\mathcal I_1 & \lesssim \mu^{-2}\lambda_1\lambda_2\beta \lambda_1^\frac3q \beta^\frac2q (\beta N_2)^{\frac12-2\eta}\lambda_2^{1-\frac3q}N_2^{\frac12+\eta}d^{-\frac12}\|w_{\mu,N}\|_{V^2_\theta}\|u_{\lambda_1,N_1}\|_{V^2_{\theta_1}}\|v_{\lambda_2,N_2}\|_{U^2_{\theta_2}}.	
\end{align*}
The summation with respect to $d\lesssim\lambda_1$ yields 
$$
\sum_{d\lesssim\lambda_1}\mathcal I_1 \lesssim \lambda_1^\frac12 \left(\frac{\lambda_1}{\lambda_2}\right)^\frac3q \left(\frac{\lambda_2}{\mu}\right)^{-2}N_2^{1-\eta}\|w_{\mu,N}\|_{V^2_\theta}\|u_{\lambda_1,N_1}\|_{U^2_{\theta_1}}\|v_{\lambda_2,N_2}\|_{U^2_{\theta_2}},	
$$
where we used the continuous embedding $U^2\subset V^2$ for $u_{\lambda_1,N_1}$. (See \cite[Proposition 2.4]{haheko}.)
This completes the proof of \eqref{main-bi-loc} in the low-modulation regime. 
\subsection{High modulation}\label{sec:high-mod}
From now on we shall consider the high-modulation region: $d\gg\min\{\mu,\lambda_1,\lambda_2\}$. In this regime we only consider low-output interaction, i.e., $\mu\ll\lambda_1\approx\lambda_2$; the Fourier multiplier $|\nabla|^{-2}$ yields good decay in the high-output interaction and hence it is much easier than the low-output case. As Section \ref{sec:low-mod}, we put $\lambda_1=\lambda_2=\lambda$. 
 Note that the angle between the Fourier support of $u_\lambda$ and $v_\lambda$ is less than $\dfrac\mu\lambda$. We put $\alpha = \dfrac\mu\lambda$. By the orthogonal decompositions by smaller cubes of size $\mu$ and conic sectors of size $\alpha$ we follow the similar approach as we have done in Section \ref{sec:low-mod}. For $d\gtrsim\lambda$, we have
\begin{align*}
\mathcal I_0 & \lesssim \sum_{\substack{ \mathtt q_1,\mathtt q_2\in \mathcal Q_\mu \\ |\mathtt q_1+\mathtt q_2|\lesssim\mu  }}\sum_{\substack{ \kappa_1,\kappa_2\in\mathcal C_\alpha \\ |\kappa_1+\kappa_2|\lesssim\alpha }}\int_{\mathbb R^{3+1}}C^\theta_d w_{\mu,N}|\nabla|^{-2}Q(P_{\mathtt q_1}R_{\kappa_1}C^{\theta_1}_{\ll d}u_{\lambda,N_1},P_{\mathtt q_2}R_{\kappa_2}C^{\theta_2}_{\ll d}v_{\lambda,N_2})\,dtdx \\
& \lesssim d^{-\frac12} \mu^{-2}\alpha \lambda^2 \|w_{\mu,N}\|_{V^2_\theta} \bigg(\sum_{\substack{ \mathtt q_1,\mathtt q_2\in \mathcal Q_\mu \\ |\mathtt q_1+\mathtt q_2|\lesssim\mu  }}\sum_{\substack{ \kappa_1,\kappa_2\in\mathcal C_\alpha \\ |\kappa_1+\kappa_2|\lesssim\alpha }}\|P_{\mathtt q_1}R_{\kappa_1}C^{\theta_1}_{\ll d}u_{\lambda,N_1}\|_{L^2_tL^\infty_x}^2\|P_{\mathtt q_2}R_{\kappa_2}C^{\theta_2}_{\ll d}v_{\lambda,N_2}\|_{L^\infty_tL^2_x}^2\bigg)^\frac12 \\
& \lesssim d^{-\frac12} \mu^{-2}\alpha \lambda^2 \mu^\frac3q\alpha^\frac2q \|w_{\mu,N}\|_{V^2_\theta}\sup_{\kappa_1\in\mathcal C_\alpha}\|R_{\kappa_1}C^{\theta_1}_{\ll d}u_{\lambda,N_1}\|_{L^2_tL^q_x} \\
& \qquad\qquad \times \bigg(\sum_{\substack{ \mathtt q_1,\mathtt q_2\in \mathcal Q_\mu \\ |\mathtt q_1+\mathtt q_2|\lesssim\mu  }}\sum_{\substack{ \kappa_1,\kappa_2\in\mathcal C_\alpha \\ |\kappa_1+\kappa_2|\lesssim\alpha }}\|P_{\mathtt q_2}R_{\kappa_2}C^{\theta_2}_{\ll d}v_{\lambda,N_2}\|_{L^\infty_tL^2_x}^2\bigg)^\frac12 \\
& \lesssim d^{-\frac12} \mu^{-2}\alpha \lambda_1\lambda_2 \mu^\frac3q\alpha^\frac2q(\alpha N_1)^{\frac12-2\eta}\lambda_1^{1-\frac3q} N_1^{\frac12+\eta} \|w_{\mu,N}\|_{V^2_\theta} \|u_{\lambda,N_1}\|_{V^2_{\theta_1}}\|v_{\lambda,N_2}\|_{V^2_{\theta_2}} \\
& \lesssim \lambda^\frac12 \left(\frac\mu\lambda\right)^{-\frac12+\frac5q+\eta}\left(\frac\lambda d\right)^\frac12 N_1^{1-\eta} \|w_{\mu,N}\|_{V^2_\theta} \|u_{\lambda,N_1}\|_{V^2_{\theta_1}}\|v_{\lambda,N_2}\|_{V^2_{\theta_2}},
\end{align*}
which gives the required bound after the summation with respect to the modulation $d;d\gtrsim\lambda$. On the other hand,
if $\mu\ll d\ll\lambda$, we see that
\begin{align*}
\mathcal I_0 & \lesssim \sum_{\substack{ \mathtt q_1,\mathtt q_2\in \mathcal Q_\mu \\ |\mathtt q_1+\mathtt q_2|\lesssim\mu  }}\sum_{\substack{ \kappa_1,\kappa_2\in\mathcal C_\alpha \\ |\kappa_1+\kappa_2|\lesssim\alpha }}\int_{\mathbb R^{3+1}}C^\theta_d w_{\mu,N}|\nabla|^{-2}Q(P_{\mathtt q_1}R_{\kappa_1}C^{\theta_1}_{\ll d}u_{\lambda,N_1},P_{\mathtt q_2}R_{\kappa_2}C^{\theta_2}_{\ll d}v_{\lambda,N_2})\,dtdx \\
& \lesssim d^{-\frac12} \mu^{-2}\alpha \lambda^2 \|w_{\mu,N}\|_{V^2_\theta} \bigg(\sum_{\substack{ \mathtt q_1,\mathtt q_2\in \mathcal Q_\mu \\ |\mathtt q_1+\mathtt q_2|\lesssim\mu  }}\sum_{\substack{ \kappa_1,\kappa_2\in\mathcal C_\alpha \\ |\kappa_1+\kappa_2|\lesssim\alpha }}\|P_{\mathtt q_1}R_{\kappa_1}C^{\theta_1}_{\ll d}u_{\lambda,N_1}\|_{L^2_tL^\infty_x}^2\|P_{\mathtt q_2}R_{\kappa_2}C^{\theta_2}_{\ll d}v_{\lambda,N_2}\|_{L^\infty_tL^2_x}^2\bigg)^\frac12 \\
& \lesssim d^{-\frac12} \mu^{-2}\alpha \lambda^2 \mu^\frac3q\alpha^\frac2q \|w_{\mu,N}\|_{V^2_\theta}\sup_{\kappa_1\in\mathcal C_\alpha}\|R_{\kappa_1}C^{\theta_1}_{\ll d}u_{\lambda,N_1}\|_{L^2_tL^q_x} \\
& \qquad\qquad \times \bigg(\sum_{\substack{ \mathtt q_1,\mathtt q_2\in \mathcal Q_\mu \\ |\mathtt q_1+\mathtt q_2|\lesssim\mu  }}\sum_{\substack{ \kappa_1,\kappa_2\in\mathcal C_\alpha \\ |\kappa_1+\kappa_2|\lesssim\alpha }}\|P_{\mathtt q_2}R_{\kappa_2}C^{\theta_2}_{\ll d}v_{\lambda,N_2}\|_{L^\infty_tL^2_x}^2\bigg)^\frac12 \\
& \lesssim d^{-\frac12} \mu^{-2}\alpha \lambda^2 \mu^\frac3q\alpha^\frac2q(\alpha N_1)^{\frac12-2\eta}\lambda^{1-\frac3q} N_1^{\frac12+\eta} \|w_{\mu,N}\|_{V^2_\theta} \|u_{\lambda_1,N_1}\|_{V^2_{\theta_1}}\|v_{\lambda_2,N_2}\|_{V^2_{\theta_2}} \\
& \lesssim \lambda^\frac12 \left(\frac\mu\lambda\right)^{-1+\frac5q+\eta}\left(\frac\mu d\right)^\frac12 N_1^{1-\eta} \|w_{\mu,N}\|_{V^2_\theta} \|u_{\lambda_1,N_1}\|_{V^2_{\theta_1}}\|v_{\lambda_2,N_2}\|_{V^2_{\theta_2}},
\end{align*}
and the summation with respect to the modulation $d;\mu\ll d\ll\lambda$ gives the desired estimate. The estimates of $\mathcal I_1$ and $\mathcal I_2$ follow by the similar way. We omit the details. This completes the proof of the main trilinear estimates \eqref{main-bi-loc}.
\section{Trilinear estimates: Proof of Theorem \ref{gwp-dirac}}
This section is devoted to the proof of Theorem \ref{gwp-dirac}. As the previous section, we define the Duhamel integral
$$
\mathfrak J^\theta[F] = \int_0^t e^{-\theta i(t-t')\langle\nabla\rangle_m}F(t')\,dt'.
$$
Then $\mathfrak J^\theta[F]$ solves the equation
$$
(-i\partial_t+\theta\langle\nabla\rangle_m)\mathfrak J^\theta[F] = F,
$$
with vanishing data at $t=0$. From now on we put $m=1$ for simplicity. We are left to prove the following trilinear estimates
\begin{align}
\|\mathfrak J^\theta [V_b*(\varphi^\dagger\phi)\psi]\|_{F^{0,1}_\theta} \lesssim \|\varphi\|_{F^{0,1}_{\theta_1}}	\|\phi\|_{F^{0,1}_{\theta_2}}\|\psi\|_{F^{0,1}_{\theta_3}}
\end{align}
which imply the global well-posedness and scattering in the $U^2$-space provided that the smallness condition for the inital data is given. The use of duality in $U^2-V^2$ gives 
\begin{align*}
\|\mathfrak J^{\theta_4}[V_b*(\varphi^\dagger\phi)\psi]\|_{F^{0,1}_{\theta_4}}^2 & \lesssim \sum_{\lambda_4,N_4\ge1}(N_4)^2\|P_{\lambda_4}H_{N_4}\mathfrak J^{\theta_4}[V_b*(\varphi^\dagger\phi)\psi]\|_{U^2_{\theta_4}}^2 \\
& \lesssim \sum_{\lambda_4,N_4\ge1}	(N_4)^2 \sup_{\|P_{\lambda_4}H_{N_4}\psi\|_{V^2_{\theta_4}}\le1}\left| \int_{\mathbb R^{1+3}} V_b*(\varphi^\dagger\phi)(P_{\lambda_4}H_{N_4}\psi)^\dagger\psi\,dtdx \right|^2. 
\end{align*}
Then dyadic decompositions and the H\"older inequality yield 
\begin{align*}
	&\|\mathfrak J^{\theta_4}[V_b*(\varphi^\dagger\phi)\psi]\|_{F^{0,1}_{\theta_4}}^2 \\
	& \lesssim \sum_{\lambda_j\ge1,j=0,1,\cdots,4}\sum_{N_j\ge1,j=0,1,\cdots,4} \\
	&\qquad \sup_{\|P_{\lambda_4}H_{N_4}\psi\|_{V^2_{\theta_4}}\le1}\left| \int_{\mathbb R^{1+3}} \langle\nabla\rangle^{-2}_bP_{\lambda_0}H_{N_0}(\varphi_{\lambda_1,N_1}^\dagger\phi_{\lambda_2,N_2})P_{\lambda_0}H_{N_0}(\psi_{\lambda_4,N_4}^\dagger\psi_{\lambda_3,N_3})\,dtdx \right|^2 \\
	& \lesssim \sum_{\lambda_j\ge1,j=0,1,\cdots,4}\sum_{N_j\ge1,j=0,1,\cdots,4} \langle\lambda_0\rangle^{-2} \\
	&\qquad \sup_{\|P_{\lambda_4}H_{N_4}\psi\|_{V^2_{\theta_4}}\le1} \|P_{\lambda_0}H_{N_0}(\varphi_{\lambda_1,N_1}^\dagger\phi_{\lambda_2,N_2})\|_{L^2_tL^2_x}^2\|P_{\lambda_0}H_{N_0}(\psi_{\lambda_4,N_4}^\dagger\psi_{\lambda_3,N_3})\|_{L^2_tL^2_x}^2.
\end{align*}
Thus our main trilinear estimates follow from the following frequency-localised $L^2$-bilinear estimates:
\begin{lem}
Let $0<\eta\ll1$ be a small positive number. For some $\frac18\le \delta\le\frac14$, we have
\begin{align}\label{main-tri-loc}
\begin{aligned}
	&\|P_{\lambda_0} (\varphi_{\lambda_1,N_1}^\dagger\phi_{\lambda_2,N_2})\|_{L^2_tL^2_x} \\
	&\quad \lesssim \lambda_0 \left(\frac{\min\{\lambda_0,\lambda_1,\lambda_2\}}{\max\{\lambda_0,\lambda_1,\lambda_2\}}\right)^\delta (\min\{N_1,N_2\})^{1-\eta}\|\varphi_{\lambda_1,N_1}\|_{U^2_{\theta_1}}\|\phi_{\lambda_2,N_2}\|_{U^2_{\theta_2}}.
\end{aligned}	
\end{align}
\end{lem}
To prove \eqref{main-tri-loc} we need to deal with the frequency interactions:
\begin{align*}
\lambda_0\ll \lambda_1\approx\lambda_2, \ \lambda_1 \ll \lambda_0\approx\lambda_2, \ 	\lambda_2 \ll \lambda_0\approx\lambda_1.
\end{align*}
Then it suffices to consider the bilinear estimates 
$$
\|P_\mu (\varphi_{\lambda,N_1}^\dagger\phi_{\lambda,N_2})\|_{L^2_tL^2_x}, \ \|P_{\lambda}(\varphi^\dagger_{\mu,N_1}\phi_{\lambda,N_2})\|_{L^2_tL^2_x}
$$
for $\mu\ll\lambda$. We first consider the first bilinear form. As the proof of Theorem \ref{gwp-wave} we apply the orthogonal decomposition of cubes of size $\mu$ and conic sectors of size $\alpha$ with $\alpha=\frac\mu\lambda$ and we use in order the H\"older inequality and the Bernstein inequality and then Lemma \ref{ang-con} and the Strichartz estimates \eqref{stri-ang} for $\varphi_{\lambda,N_1}$ 
\begin{align*}
\|P_\mu (\varphi^\dagger_{\lambda,N_1}\phi_{\lambda,N_2})\|_{L^2_tL^2_x} & \lesssim \bigg(\sum_{\substack{\mathtt q_1,\mathtt q_2\in\mathcal Q_\mu \\ |\mathtt q_1-\mathtt q_2|\lesssim\mu}}\sum_{\substack{\kappa_1,\kappa_2\in\mathcal C_\alpha \\ |\kappa_1-\kappa_2|\lesssim\alpha}}\|	P_\mu (P_{\mathtt q_1}R_{\kappa_1}\varphi^\dagger_{\lambda,N_1}P_{\mathtt q_2}R_{\kappa_2}\phi_{\lambda,N_2})\|_{L^2_tL^2_x}^2 \bigg)^\frac12 \\
& \lesssim \bigg(\sum_{\substack{\mathtt q_1,\mathtt q_2\in\mathcal Q_\mu \\ |\mathtt q_1-\mathtt q_2|\lesssim\mu}}\sum_{\substack{\kappa_1,\kappa_2\in\mathcal C_\alpha \\ |\kappa_1-\kappa_2|\lesssim\alpha}}\|	 P_{\mathtt q_1}R_{\kappa_1}\varphi_{\lambda,N_1}\|^2_{L^2_tL^\infty_x}\|P_{\mathtt q_2}R_{\kappa_2}\phi_{\lambda,N_2}\|_{L^\infty_tL^2_x}^2 \bigg)^\frac12 \\
& \lesssim \mu^\frac3q \left(\frac\mu\lambda\right)^\frac2q\sup_{\kappa_1}\|R_{\kappa_1}\varphi_{\lambda,N_1}\|_{L^2_tL^q_x} \\
& \qquad\qquad\times \bigg(\sum_{\substack{\mathtt q_1,\mathtt q_2\in\mathcal Q_\mu \\ |\mathtt q_1-\mathtt q_2|\lesssim\mu}}\sum_{\substack{\kappa_1,\kappa_2\in\mathcal C_\alpha \\ |\kappa_1-\kappa_2|\lesssim\alpha}}\|P_{\mathtt q_2}R_{\kappa_2}\phi_{\lambda,N_2}\|_{L^\infty_tL^2_x}^2 \bigg)^\frac12 \\
& \lesssim \mu^\frac3q \left(\frac\mu\lambda\right)^\frac2q \left(\frac\mu\lambda N_1\right)^{\frac12-2\eta}\lambda^{1-\frac3q}N_1^{\frac12+\eta}\|\varphi_{\lambda,N_1}\|_{U^2_{\theta_1}}\|\phi_{\lambda,N_2}\|_{U^2_{\theta_2}} \\
& \lesssim \mu \left(\frac\mu\lambda\right)^{\frac5q-\frac12-\eta}N_1^{1-\eta} \|\varphi_{\lambda,N_1}\|_{U^2_{\theta_1}}\|\phi_{\lambda,N_2}\|_{U^2_{\theta_2}}.
\end{align*}
If $N_2\ll N_1$, then we interchange the role of $\varphi$ and $\phi$.
 For the second bilinear form, we are only concerned with the case $N_1\gg N_2$. We make the use of $L^2$-duality and then orthgonal decompositions of cubes and angular sectors of size $c>0$ where $c$ is a small constant and we have
\begin{align*}
\|P_{\lambda}(\varphi_{\mu,N_1}^\dagger\phi_{\lambda,N_2})\|_{L^2_tL^2_x} & \lesssim \sup_{\|\psi\|_{L^2_tL^2_x}\lesssim1}\sum_{\substack{q,q_2\in\mathcal Q_\mu \\ |q+\theta_2\mathtt q_2|\lesssim\mu}}\sum_{\substack{\kappa,\kappa_1,\kappa_2\in\mathcal C_c \\ |\kappa_1-\kappa_2|,|\kappa+\theta_2\kappa_2|\lesssim c}}\\
& \qquad\qquad \int_{\mathbb R^{1+3}}P_qR_\kappa\psi R_{\kappa_1}\varphi^\dagger_{\mu,N_1}P_{\mathtt q_2}R_{\kappa_2}\phi_{\lambda,N_2}\,dtdx	 \\
& \lesssim \sup_{\|\psi\|_{L^2_tL^2_x}\lesssim1}\bigg(\sum_{\substack{q,q_2\in\mathcal Q_\mu \\ |q+\theta_2\mathtt q_2|\lesssim\mu}}\sum_{\substack{\kappa,\kappa_1,\kappa_2\in\mathcal C_c \\ |\kappa_1-\kappa_2|,|\kappa+\theta_2\kappa_2|\lesssim c}} \\
&\qquad\qquad \|P_qR_\kappa\psi \|_{L^2_tL^2_x}^2\|R_{\kappa_1}\varphi_{\mu,N_1}P_{\mathtt q_2}R_{\kappa_2}\phi_{\lambda,N_2}\|_{L^2_tL^2_x}^2\bigg)^\frac12 \\
& \lesssim \sup_{\|\psi\|_{L^2_tL^2_x}\lesssim1}\bigg(\sum_{\substack{q,q_2\in\mathcal Q_\mu \\ |q+\theta_2\mathtt q_2|\lesssim\mu}}\sum_{\substack{\kappa,\kappa_1,\kappa_2\in\mathcal C_c \\ |\kappa_1-\kappa_2|,|\kappa+\theta_2\kappa_2|\lesssim c}} \\
& \qquad\qquad\qquad\qquad \|P_qR_\kappa\psi \|_{L^2_tL^2_x}^2\|R_{\kappa_1}\varphi_{\mu,N_1}\|_{L^\infty_tL^2_x}^2\|P_{\mathtt q_2}R_{\kappa_2}\phi_{\lambda,N_2}\|_{L^2_tL^\infty_x}^2\bigg)^\frac12 \\
& \lesssim \mu^\frac3q c^\frac2q\sup_{\kappa_2}\|R_{\kappa_2}\phi_{\lambda,N_2}\|_{L^2_tL^q_x} \|\varphi_{\mu,N_1}\|_{V^2_{\theta_1}} \\
& \lesssim \mu^\frac3q c^\frac2q(cN_2)^{\frac12-2\eta}\lambda^{1-\frac3q}N_2^{\frac12+\eta}\|\varphi_{\mu,N_1}\|_{V^2_{\theta_1}}\|\phi_{\lambda,N_2}\|_{U^2_{\theta_2}} \\
& \lesssim \lambda \left(\frac\mu\lambda\right)^\frac3q N_2^{1-\eta} \|\varphi_{\mu,N_1}\|_{u^2_{\theta_1}}\|\phi_{\lambda,N_2}\|_{U^2_{\theta_2}},
\end{align*}
where we used the continuous embedding $U^2\subset V^2$. This completes the proof of the main bilinear estimates \eqref{main-tri-loc}. 

%%%%%%%%%%%%%%%%%%%%%%%%%%%%%%%%%%%%%%%%%%%%%%%%%%%%%%%%%%%%%%%%%%%%%%%%%%%%%%%%%%%%%%%%%%%%%%%%%%%%%%%%%%%%%%%%%%%%%%%%%%%%%%%%%%%%%%%%%%%%%
\section*{Acknowledgements}
Most of all, the author would like to express his gratitude to Doctor Lee, Kiyeon, and Professor Cho, Yonggeun for helpful discussion and generous criticism. This work was supported by the National Research Foundation of Korea (NRF) grant funded by the Korea government (MSIT) (NRF-2020R1A2C4002615). 
%%%%%%%%%%%%%%%%%%%%%%%%%%%%%%%%%%%%%%%%%%%%%%%%%%%%%%%%%%%%%%%%%%%%%%%%%%%%%%%%%%%%%%%%%%%%%%%%%%%%%%%%%%%%%%%%%%%%%%%%%%%%%%%%%%%%%%%%%%%%%
%%%%%%%%%%%%%%%%%%%%%%%%%%%%%%%%%%%%%%%%%%%%%%%%%%%%%%%%%%%%%%%%%%%%%%%%%%%%%%%%%%%%%%%%%%%%%%%%%%%%%%%%%%%%%%%%%%%%%%%%%%%%%%%%%%%%%%%%%%%%%

%%%%%%%%%%%%%%%%%%%%%%%%%%%%%%%%%%%%%%%%%%%%%%%%%%%%%%%%%%%%%%%%%%%%%%%%%%%%%%%%%%%%%%%%%%%%%%%%%%%%%%%%%%%%%%%%%%%%%%%%%%%%%%%%%%%%%%%%%%%%%%%%%%%%%%%%%%%%%%%%%%%%%%%%%%%%%%%%%%%%%%%%%%%%%%%%%%%%%%%%%%%%%%%%%%%%%%%%%%%%%%%%%%%%%%%%%%%%%%%%%%%%%%%%%%%%%%%%%%%%%%%%%%%%%%%%%%%%%%%%%%%%%%%%%%%%%%%%%%%%%%%%%%%%


\begin{thebibliography}{00}

\bibitem{danfos} P. D' Ancona, D. Foschi, S. Selberg, {\it Null structure and almost optimal local regularity for the Dirac-Klein-Gordon system}, J. Eur. Math. Soc. (JEMS) \textbf{9}:4, (2007), 877--899. 

\bibitem{dasfos1}  P. D' Ancona, D. Foschi, S. Selberg, {\it Null structure and almost optimal local well-posedness of the Maxwell-Dirac system}, Amer. J. Math. Vol. 132, No. 3, (2010), 771--839.

\bibitem{behe} I. Bejenaru and S. Herr, {\it On global well-posedness and scattering for the massive Dirac-Klein-Gordon system}, Journal of European Mathematics Society, \textbf{19}:8, (2017): 2445--2467.

\bibitem{bjor} J.D. Bjorken, S. D. Drell, {\it Relativistic quantum mechanichs}, McGraw-Hill, New York, 1964.

 \bibitem{candyherr} T. Candy and S. Herr, {\it Transference of bilinear restriction estimates to quadratic variation norms and the Dirac-Klein-Gordon system}, Analysis and PDE 11, no. 5, (2018): 1171--1240.

\bibitem{candyherr1} T. Candy and S. Herr, {\it Conditional large initial data scattering results for the Dirac-Klein-Gordon system}, Forum of Mathematics, Sigma, (2018), Vol. 6, 55 pp.

\bibitem{cholee} Y. Cho and K. Lee, {\it Small data scattering of Dirac equations with Yukawa type potentials in $L^2_x(\mathbb R^2)$}, Differential Integral equations, \textbf{34}, (2021), 425--436.

\bibitem{choslee} Y. Cho, S. Lee, {\it Strichartz estimates in spherical coordinates}, Indiana University Mathematics Journal 62, no. 3, (2013): 991--1020.

\bibitem{chohlee} Y. Cho, S. Hong, and K. Lee, {\it Scattering and non-scattering of the Hartree-type nonlinear Dirac system at critical regularity}, arXiv:2106.01633.

\bibitem{choozlee} Y. Cho, K. Lee, and T. Ozawa, {\it Small data scattering of 2d Hartree type Dirac equations}, J. Math. Anal. Appl. \textbf{506}, (2022), 125549.

\bibitem{chooz}  Y. Cho, T. Ozawa, {\it On the semirelativistic Hartree-type equation}, SIAM J. Math. Anal. \textbf{38} (2006), 1060--1074.

\bibitem{gaoh} G. Critstian, S.-J. Oh, {\it Global well-posedness of high dimensional Maxwell-Dirac for small critical data}, Mem. Amer. Math. Soc. \textbf{264}, No. 1279. (2020)

%\bibitem{choozxia} Y. Cho, T. Ozawa, and S. Xia, {\it Remarks on some dispersive estimates}, Communications in Pure and Applied Analysis, \textbf{10}, (2011), 1121--1128.

\bibitem{fosklai} D. Foschi, S. Klainerman, {\it Bilinear space-time estimates for homogeneous wave equations}, Annales Scientifiques de l' \'Ecole Normale Sup\'erieure 23, No. 2, (2000): 211--74.

%\bibitem{geosha} V. Georgiev and B. Shakarov, {\it Global large data solutions for 2D Dirac equation with Hartree type interaction}, International Mathematics Research Notices, (2021) 

\bibitem{haheko} M. Hadac, S. Herr, and H. Koch, {\it Well-posedness and scattering for the KP-II equation in a critical space}, Inst. H.Poincare Anal. Non lineaire, \textbf{26}, (2009), 917--941.

\bibitem{herrlenz} S. Herr, E. Lenzmann, {\it The Boson star equation with initial data of low regularity}, Nonlinear Analysis, Vol. 97, (2014): 125--137.

\bibitem{herrtes} S. Herr, A. Tesfahun, {\it Small data scattering for semi-relativistic equations with hartree type nonlinearity}, Journal of Differential Equations, Vol. 259, (2015): 5510--5532.

\bibitem{huhoh} H. Huh, S.-J. Oh, {\it Low regularity solutions to the Chern-Simons-Dirac and the Chern-Simons-Higgs equations in the Lorenz gauge}, Communications in Partial Differential Equations, \textbf{41} (3) (2016), 375--397.

\bibitem{keeltao} M. Keel, T. Tao, {\it Endpoint Strichartz estimates}, American Journal of Mathematics, \textbf{120}, No. 5, (1998), 955--980.

\bibitem{klaima} S. Klainerman, M. Machedon, {\it Space-time estimates for null forms and the local existence theorem}, Communications in Pure and Applied Mathematics, Vol. 46, No. 9, (1993), 1221--1268.

\bibitem{klaima1} S. Klainerman, M. Machedon, {\it Smoothing estimates for null forms and applications}, Duke Mathematics Journal, \textbf{81}, (1995), 99--133.

\bibitem{klaitataru} S. Klainerman, D. Tataru, {\it On the optimal local regularity for Yang-Mills equations in $\mathbb R^{4+1}$}, Journal of American Mathematical Society 12, (1999): 93--116.

\bibitem{kochtavi} H. Koch, D. Tataru, M. Visan, {\it Dispersive equations and nonlinear waves}, Basel: Birkh\"auser/Springer, 2014.

\bibitem{kriegersterbenztataru} J. Krieger, J. Sterbenz, D. Tataru, {\it Global well-posedness for the Maxwell-Klein-Gordon equation in $4+1$ dimensions: small energy}, Duke Mathematical Journal 164, no. 6, (2015): 973--1040.

\bibitem{kriegertataru} J. Krieger, D. Tataru, {\it Global well-posedness for the Yang-Mills equation in $4+1$ dimensions. Small energy}, Annals of Mathematics, \textbf{185}, (2017), 831--893.

\bibitem{sleevar} S. Lee, A. Vargas, {\it Sharp null form estimates for the wave equation}, American Journal of Mathematics, \textbf{130}, No. 5, (2008), 1279--1326.

\bibitem{masterbenz} M. Machedon, J. Sterbenz, {\it Almost optimal local well-posedness for the $(3+1)-$dimensional Maxwell-Klein-Gordon equations}, Jounal of American Mathematical Society 17, (2004):297--359.

\bibitem{ohtataru} S.-J. Oh, D. Tataru, {\it Local well-posedness of the $(4+1)$-dimensional Maxwell-Klein-Gordon equation at energy regularity}, Annals of PDE 2, No. 1, Art. 2, (2016), 70pp.

\bibitem{ohtataru1} S.-J. Oh, D. Tataru, {\it Global well-posedness and scattering of the $(4+1)$-dimensional Maxwell-Klein-Gordon-equation}, Invent. Math. 205. No.3, (2016), 781--877.

\bibitem{ohtataru2} S.-J. Oh, D. Tataru, {\it Energy dispersed solutions for the $(4+1)$-dimensional Maxwell-Klein-Gordon equation}, Amer. J. Mathe. 140, No. 1, (2018), 1--82.

\bibitem{ster} J. Sterbenz, {\it Angular regularity and Strichartz estimates for the wave equation}, Int. Math. Res. Not. 2005:4 (2005), 187--231.

\bibitem{sterbenz2} J. Sterbenz, {\it Global regularity for general non-linear wave equations I\!I. $(4+1)-$dimensional Yang-Mills equations in the Lorenz gauge }, American Journal of Mathematics, Vol. 129, No. 3, (2007), 611--664.


\bibitem{stri} R. S. Strichartz, {\it Restrictions of Fourier transforms to quadratic surfaces and decay of solutions of wave equations}, Duke Mathematics Journal, \textbf{44}, No. 3, (1977), 705--714.

\bibitem{tao} T. Tao, {\it Low regularity semi-linear wave equations}, Communications in Partial Differential Equations 24, No. 3-4, (1999): 599--629.

\bibitem{tao1} T. Tao, {\it Global regularity of wave maps, I: small critical Sobolev norm in high dimension}, International Mathematics Research Notices 2001, No. 6, (2001): 299--328.

\bibitem{tao2} T. Tao, {\it Global regularity of wave maps I\!I. small energy in two dimensions}, Communication in Mathematical Physics 224, No. 2, (2001): 443--544.

\bibitem{tataru} D. Tataru, {\it On the equation $\Box u=|\nabla u|^2$ on $5+1$ dimensions}, Mathematical Research Letters 6, No. 5-6, (1999): 469--485.
\bibitem{tataru1} D. Tataru, {\it On global existence and scattering for the wave map equations}, American Journal of Mathematics 123, No. 1, (2001): 37--77.  

\bibitem{tes} A. Tesfahun, {\it Long-time behavior of solutions to cubic Dirac equation with Hartree type nonlinearity in $\mathbb R^{1+2}$}, International Mathematics Research Notices 2020: 19, (2020): 6489--6538.

\bibitem{tes1} A. Tesfahun, {\it Small data scattering for cubic Dirac equation with Hartree type nonlinearity in $\mathbb R^{1+3}$}, SIAM Journal of Mathematical Analysis, Vol. 52, No. 3, (2020): 2969--3003.

\bibitem{cyang} C. Yang, {\it Scattering results for Dirac Hartree-type equations with small initial data}, Communication in Pure and Applied Analysis, Vol. 18, No. 4, (2019), 1711--1734.

\bibitem{wang} X. Wang, {\it On global existence of 3D charge critical Dirac-Klein-Gordon system}, International Mathematics Reserach Notices 2015, no.21, (2015): 10801--10846.


\end{thebibliography}
\end{document}